\documentclass[a4paper,12pt]{article}

\usepackage{amsmath,amssymb,amsthm,amscd,cases}
{}
\numberwithin{equation}{section}

\newtheorem{thm}{Theorem}[section]

\newcommand{\R}{\mathbb{R}}

\renewcommand{\a}{\alpha}

\renewcommand{\d}{\delta}
\newcommand{\p}{\partial}

\newcommand{\bc}{\begin{cor}}
\newcommand{\ec}{\end{cor}}
\newcommand{\bl}{\begin{lem}}
\newcommand{\el}{\end{lem}}
\newcommand{\bp}{\begin{prop}}
\newcommand{\ep}{\end{prop}}
\newcommand{\bt}{\begin{thm}}
\newcommand{\et}{\end{thm}}
\newcommand{\bal}{\begin{array}{ll}}
\newcommand{\ba}{\begin{array}}
\newcommand{\bac}{\begin{array}{ccc}}
\newcommand{\ea}{\end{array}}

\newcommand{\be}{\begin{equation}}
\newcommand{\ee}{\end{equation}}

\newtheorem{lem}[thm]{Lemma}
\newtheorem{cor}[thm]{Corollary}
\newtheorem{prop}[thm]{Proposition}

\begin{document}

\title{\Large\bf {Asymptotic Behavior of Solutions to a Model System of a
Radiating Gas}}

\author{Yongqin Liu{\footnote{email:~yqliu2@yahoo.com.cn}}\\
 \small\emph{Faculty of Mathematics, Kyushu University}\\
 \small\emph{Fukuoka 819-0395, Japan} \\[3mm]
Shuichi Kawashima{\footnote{email:~kawashim@math.kyushu-u.ac.jp}}\\
 \small\emph{Faculty of Mathematics, Kyushu University}\\
 \small\emph{Fukuoka 819-0395, Japan}
 }

\date{}

\maketitle

 \begin{abstract}

In this paper we focus on the initial value problem for a hyperbolic-elliptic 
coupled system of a radiating gas in multi-dimensional space. 
By using a time-weighted energy method, we obtain the global existence and 
optimal decay estimates of solutions. Moreover, we show that the solution 
is asymptotic to the linear diffusion wave which is given in terms of the 
heat kernel. \\

\noindent 
{\bf keywords}: radiating gas; initial value problem; asymptotic behavior. \\


\end{abstract}

\section{Introduction}

In this paper we consider the initial value problem: 
\be\label{1a}
\left\{\begin{array}{ll}
u_{t}+\sum\limits_{j=1}\limits^{n}f_j(u)_{x_j}
   +{\rm div}q=0,\ &x\in\R^n,\ \ t>0, \\[2mm]
-\nabla{\rm div}q+q+\nabla u=0,\ &x\in\R^n,\ \ t>0,
\end{array}
\right. 
\ee
with  the initial data 
\be\label{1a1}
u(x,0)=u_0(x),\quad x\in\R^n.
\ee
Here $f_j(u),\ j=1,\cdots,n$, are smooth functions of $u$ satisfying 
$f_j(u)=O(u^2)$ for $u\to 0$, and $u=u(x,t)$ and 
$q=(q_1,\cdots,q_n)(x,t)$ are unknown functions of 
$x=(x_1,\cdots,x_n)\in\R^n\ (n\geq1)$ and $t>0$. 
Typically, $u$ and $q$ represent the velocity and radiating heat flux of 
the gas, respectively. 

The system \eqref{1a} is a simplified version of a radiating gas model 
in $n$-dimensional space.  More precisely, 
in a certain physical situation, the system \eqref{1a} gives a good
approximation to the following system of a radiating gas, 
that is a quite general model for compressible gas dynamics where 
the heat radiative transfer phenomena are taken into account: 
\be\label{1b}
\left\{\bal 
&\rho_t+{\rm div}(\rho u)=0,\\[1mm]
&(\rho u)_t+{\rm div}(\rho u\otimes u+pI)=0,\\[1mm]
&\{
\rho(e+{{|u|^2}\over2})\}_t
   +{\rm div}\{\rho u(e+{{|u|^2}\over2})+pu+q\}=0,\\[1mm]
&-\nabla{\rm div}q+a_1q+a_2\nabla\theta^4=0, 
\ea
\right. 
\ee
where $\rho,\ u,\ p,\ e$ and $\theta$ are respectively the mass density, 
velocity, pressure, internal energy and absolute temperature of the 
gas, while $q$ is the radiative heat flux, and $a_1$ and $a_2$ are 
given positive constants depending on the gas itself. 
The first three equations form the usual Euler system, which 
describes the inviscid flow of a compressible fluid, and express 
the conservation of mass, momentum and energy, respectively. 
We refer to the book of Courant and Friedrichs \cite{CF} for a detailed
derivation of several models in compressible gas dynamics. 
On the other hand, the physical motivation of the fourth equation, 
which takes into account the heat radiation phenomena, 
is given in \cite{VK,Di1}. 
The simplified model \eqref{1a} was first recovered by Hamer \cite{Ha}, 
and the reduction of the full system \eqref{1b} to \eqref{1a} was 
given in \cite{VK,Ha,KT, GZ}. 

There are many works on the study of the hyperbolic-elliptic coupled 
system for one-dimensional radiating gas. The earlier paper with application to this kind of systems is \cite{ST}, where  Schochet and  Tadmor studied the regularized Chapman-Enskog expansion for scalar conservation laws.
We refer to \cite{Ha,KN2,KN,KN1, LCG,LCG1,LMS,LMS1,LMNPZ,LT,NPZ,Ni} for shock waves, 
\cite{DFZ,KNN,IK,KNN2} for nonlinear diffusion waves, 
\cite{KT} for rarefaction waves, 
\cite{KN3,La,LM} for a singular limit and relaxation limit, 
and \cite{DL,Se,Se1} for $L^1$ stability results, . 
In the  multi-dimensional case, 
Di Francesco in \cite{Di} obtained the global-in-time existence and 
uniqueness of weak entropy solutions to the system \eqref{1a} 
and analyzed the relaxation limits. 
In \cite{WW}, recently Wang and Wang studied the initial value problem for the
system \eqref{1a} in multi-dimensions and obtained the pointwise estimates of classical solutions  by using the method of Green function combined with some energy estimates.
Also, the stability of planar rarefaction waves was discussed in 
\cite{GZ,GRZ}. 
Very recently, Ruan and Zhu \cite{RZ} investigated the asymptotic 
decay rates toward the linear diffusion wave and also the rarefaction wave in $\mathbb{R}^n$, and proved the asymptotic relation of our Proposition \ref{NL} with $k=0$ and the norm $L^2$ instead of $H^{s-1-k}$ in the case $2\leq n\leq 7$ by using the energy method and the semigroup argument. 

In this paper we investigate the decay rate not only to the same linear diffusion wave as in \cite{RZ} which can be seen from Proposition \ref{NL}, but also to another diffusion wave which is given in terms of  the heat kernel as shown in Theorem \ref{22}
for the initial value problem \eqref{1a}, \eqref{1a1} in $\mathbb{R}^n$ 
for $n\geq 2$ by applying the time-weighted energy method together with 
the semigroup argument , which removes the restriction $n\leq 7$ assumed in 
\cite{RZ} and also improves their results. 
To this end, we first transform the system \eqref{1a} into the following equivalent 
decoupled system \eqref{3a} which makes the derivation of our energy 
estimates easier, but is not essential for obtaining our main results,
\be\label{3a}
\left\{\bal u_t-\Delta u_t-\Delta u
+(1-\Delta)\sum\limits^n\limits_{j=1}f_j(u)_{x_j}=0,\ 
&x\in\R^n,\ t>0,\\
q=-(1-\Delta)^{-1}\nabla u,\ &x\in\R^n,\ t>0.
\ea\right.
\ee 
 Then, by applying the time-weighted energy method to 
the decoupled system \eqref{3a}, we derive the optimal decay estimates 
of solutions for all $n\geq 1$. 
Finally, using the semigroup argument, we show that the solution is 
asymptotic to the linear diffusion wave as $t\to+\infty$, provided that 
$n\geq 2$. Our linear diffusion wave is given explicitly in terms of 
the heat kernel. 

The contents of the paper are as follows. 
In Section 2 we give full statements of our main theorems.  
Section 3 gives the proof of the results on the global existence 
and decay estimates of solutions. 
The last section gives the proof of the theorem on the asymptotic 
convergence to the linear diffusion wave. 

Before closing this section, we give some notations to be used 
below. Let $\mathcal{F}[f]$ denote the Fourier transform of $f$ 
defined by 
$$
\mathcal{F}[f](\xi)=\hat{f}(\xi):=
\int_{\R^n}e^{-ix\cdot \xi}f(x)dx,
$$
and we denote its inverse transform by $\mathcal{F}^{-1}.$ 

$L^p=L^p(\R^n)\ (1\leq p\leq+\infty)$ is the usual Lebesgue space with 
the norm $\|\cdot\|_{L^p}$. 
For $\gamma\in\R$, let $L^1_{\gamma}=L^1_{\gamma}(\R^n)$ denote the
weighted $L^1$ space with the norm 
$$
\|f\|_{L^1_{\gamma}}:=\int_{\R^n}(1+|x|)^{\gamma}|f(x)|dx.
$$
Let $s$ be a nonnegative integer. 
Then $H^s=H^s(\R^n)$ denotes the Sobolev space of $L^2$ functions, 
equipped with the norm 
$$
\|f\|_{H^s}:=\Big(\sum\limits_{k=0}\limits^{s}
\|\partial_x^kf\|_{L^2}^2\Big)^{1\over2}.
$$
Here, for a nonnegative integer $k$, $\partial_x^k$ denotes the totality 
of all the $k$-th order derivatives with respect to 
$x\in{\mathbb R}^n$. 
Also, $C^k(I; H^s(\R^n))$ denotes the space of $k$-times 
continuously differentiable functions on the interval $I$ 
with values in the Sobolev space $H^s=H^s(\R^n)$. 

Finally, in this paper, we denote every positive constant by the same 
symbol $C$ or $c$ without confusion. $[\,\cdot\,]$ is the Gauss' symbol. 

\section{Main theorems}

Our first theorem is on the global existence and uniform energy estimate 
of solutions to the problem \eqref{1a}, \eqref{1a1}. 

\bt\label{21} 
Let $n\geq1$ and let $s\geq[{n\over2}]+2$ be an integer. 
Assume that $u_0\in H^{s}(\R^n)$ and put $E_0:=\|u_0\|_{H^{s}}.$ 
Then there exists a small positive constant $\d_0$ such that 
if $E_0\leq\d_0,$ then the problem \eqref{1a}, \eqref{1a1} 
has a unique global solution $(u,q)(x,t)$ with 
\begin{eqnarray*}
&u\in C([0,+\infty); H^{s}(\R^n)), \quad
   \partial_x u\in L^2([0,+\infty); H^{s-1}(\R^n)),\\[1mm]
&q\in C([0,+\infty); H^{s+1}(\R^n))\cap L^2([0,+\infty); H^{s+1}(\R^n)).
\end{eqnarray*}
The solution verifies the uniform energy estimate 
\begin{equation}\label{2.1}
\|u(t)\|_{H^s}^2+\|q(t)\|_{H^{s+1}}^2
+\int_0^t\|\partial_xu(\tau)\|_{H^{s-1}}^2
+\|q(\tau)\|_{H^{s+1}}^2d\tau
\leq CE_0^2.
\end{equation}
\et

When $n\geq 2$, the solution obtained in Theorem \ref{21} verifies 
the following decay estimates. 

\bc\label{21a}
Let $n\geq2$ and let $s\geq[{n\over2}]+2$ be an integer. 
If $E_0=\|u_0\|_{H^{s}}$ is suitably small, then the global solution 
obtained in Theorem \ref{21} verifies the decay estimates 
\begin{equation}\label{2.2}
\|\p^k_xu(t)\|_{H^{s-k}}\leq CE_0(1+t)^{-{k\over2}}
\end{equation}
for $k$ with $0\leq k\leq s$ and 
\begin{equation}\label{2.3}
\|\p^k_xq(t)\|_{H^{s+1-k}}\leq CE_0(1+t)^{-{k+1\over2}}
\end{equation}
for $k$ with $0\leq k\leq s-1.$
\ec

The next theorem is on the optimal decay estimates of solutions 
for initial data in $H^{s}(\R^n)\cap L^1(\R^n)$. 

\bt\label{21b}
Let $n\geq1$, and let $s\geq3$ for $n=1$ and 
$s\geq[{n\over2}]+2$ for $n\geq2$. 
Assume that $u_0\in H^{s}(\R^n)\cap L^1(\R^n)$ and put 
$E_1:=\|u_0\|_{H^{s}}+\|u_0\|_{L^1}.$ 
Then there exists a small positive constant $\d_1$ such that 
if $E_1\leq\d_1,$ then the global solution 
obtained in Theorem \ref{21} satisfies the decay estimates 
\begin{equation}\label{2.4}
\|\p^k_xu(t)\|_{H^{s-k}}\leq CE_1(1+t)^{-{n\over4}-{k\over2}}
\end{equation}
for $k$ with $0\leq k\leq s-1$ and 
\begin{equation}\label{2.5}
\|\p^k_xq(t)\|_{H^{s+1-k}}\leq CE_1(1+t)^{-{n\over4}-{k+1\over2}}
\end{equation}
for $k$ with $0\leq k\leq s-2.$
\et

\noindent
{\bf Remark 1.}\  
To obtain the optimal decay estimates of solutions stated in 
Corollary \ref{21a} and Theorem \ref{21b}, we will use the following decay estimate for $L^{\infty}$ norm of 
the derivative $\p_xu$ as shown by $N(T)$ in \eqref{N}:
$$
\|\p_xu(t)\|_{L^{\infty}}\leq C(1+t)^{-1}.
$$

Our final result is concerning the asymptotic profile of the global 
solution obtained in Theorem \ref{21b} for $n\geq 2$. 
First we show that for $n\geq 2$, the solution to the problem 
\eqref{1a}, \eqref{1a1} can be approximated by the solution to 
the corresponding linear problem,
\be\label{32a}
\left\{\bal 
&\bar{u}_t-\Delta \bar{u}_t-\Delta\bar{u}=0, \quad
\bar{u}(x,0)=u_0(x),
\\[2mm]
&\bar{q}=-(1-\Delta)^{-1}\nabla\bar{u}.
\ea \right. 
\ee 
Then we prove that the solution to this linear problem can be further 
approximated by the following simpler problem \eqref{HE} based on the linear 
heat equation $u_t-\Delta u=0$,
\begin{equation}\label{HE}
\left\{\begin{array}{ll}
&\tilde{u}_t-\Delta \tilde{u}=0, \quad 
   \tilde{u}(x,0)=u_0(x), \\[2mm]
&\tilde{q}=-(1-\Delta)^{-1}\nabla\tilde{u}.
\end{array}\right.
\end{equation}
Since the solution to the linear heat equation is asymptotic to the 
heat kernel 
\begin{equation}\label{FS}
G_0(x,t)=\mathcal{F}^{-1}[e^{-|\xi|^2t}](x)
=(4\pi t)^{-{n\over2}}e^{-{|x|^2\over4t}},
\end{equation}
we thus conclude that the asymptotic profile of our global solution is 
given by the following linear diffusion wave $(u^*,q^*)(x,t)$: 
\begin{equation}\label{LDW}
u^*(x,t)=MG_0(x,t+1), \quad
q^*(x,t)=-\nabla u^*(x,t), 
\end{equation}
where $M=\int_{{\mathbb R}^n}u_0(x)dx$ denotes the "mass". 
%
%
The result is precisely stated as follows. 

\bt\label{22}  
Let $n\geq 2$ and $s\geq [{n\over2}]+2$. 
Assume that $u_0\in H^{s}(\R^n)\cap L^1_1(\R^n)$ and put 
$E_2:=\|u_0\|_{H^{s}}+\|u_0\|_{L^1_1}$. 
Then the global solution $(u,q)$ to the problem \eqref{1a}, \eqref{1a1}, 
which is constructed in Theorem \ref{21b}, is asymptotic to the 
linear diffusion wave $(u^*,q^*)$ in \eqref{LDW} as $t\to+\infty$: 
\begin{equation}\label{asymp-1}
\|\p^k_x(u-u^*)(t)\|_{H^{s-1-k}}\leq
   CE_2\rho(t)(1+t)^{-{n\over4}-{{k+1}\over2}}
\end{equation}
for $0\leq k\leq s-1$ and 
\begin{equation}\label{asymp-2}
\|\p^k_x(q-q^*)(t)\|_{H^{s-k}}\leq
   CE_2\rho(t)(1+t)^{-{n\over4}-{{k+2}\over2}}
\end{equation}
for $0\leq k\leq s-2$, where $\rho(t)=\ln(1+t)$ for $n=2$ and 
$\rho(t)=1$ for $n\geq 3$. 
\et

%
\section{Global existence and decay estimates}
\subsection{Global existence of solutions}

This subsection is devoted to the proof of the global existence result 
stated in Theorem \ref{21}. Since a local existence result can be 
obtained by the standard method based on the successive approximation 
sequence, we omit its details and only derive the desired a priori 
estimates of solutions. 

First we give a lemma which will be used in the derivation of our 
energy estimates. 

\bl\label{31}
Let $1\leq p,\,q,\,r\leq+\infty$ and 
${1\over p}={1\over q}+{1\over r}.$ Then the following estimates hold: 
$$
\|\p^k_x(uv)\|_{L^p}\leq
C(\|u\|_{L^q}\|\p^k_xv\|_{L^r}+\|v\|_{L^q}\|\p^k_xu\|_{L^r}), 
\quad k\geq 0.
$$
$$
\|[\p^k_x,u]\p_xv\|_{L^p}\leq
C(\|\p_xu\|_{L^q}\|\p^k_xv\|_{L^r}+\|\p_xv\|_{L^q}\|\p^k_xu\|_{L^r}), 
\quad k\geq 1.
$$
\el

The proof of this lemma can be found in \cite{LK}. 

The next lemma shows the equivalence of the system \eqref{1a} and \eqref{3a}. 

\bl\label{310}
The system \eqref{1a} is equivalent to the  decoupled system \eqref{3a}.
\el
\begin{proof}
First we show that \eqref{1a} implies \eqref{3a}. 
It follows from $\eqref{1a}_2$ that $q=\nabla\phi$ with 
$\phi={\rm div}q-u$. Therefore we see that 
$\nabla{\rm div}q=\nabla\Delta\phi=\Delta\nabla\phi=\Delta q$. 
Consequently, $\eqref{1a}_2$ becomes  
$(1-\Delta)q+\nabla u=0$, which gives $\eqref{3a}_2$. 
Also, applying ${\rm div}$ to $\eqref{3a}_2$, we have 
\begin{equation}\label{3.2}
{\rm div}q=-(1-\Delta)^{-1}\Delta u.
\end{equation}
We substitute \eqref{3.2} into $\eqref{1a}_1$ and then apply 
$1-\Delta$ to the resulting equation. This yields $\eqref{3a}_1$. 

Next we derive \eqref{1a} from \eqref{3a}. 
We have \eqref{3.2} from $\eqref{3a}_2$. Apply $(1-\Delta)^{-1}$ to 
$\eqref{3a}_1$ and substitute \eqref{3.2} to the result. This yields 
$\eqref{1a}_1$. 
Also, we rewrite \eqref{3.2} as ${\rm div}q=u-(1-\Delta)^{-1}u$ and 
apply $\nabla$, obtaining 
$$
\nabla{\rm div}q=\nabla u-(1-\Delta)^{-1}\nabla u=\nabla u+q,
$$
where we have used $\eqref{3a}_2$. This gives $\eqref{1a}_2$. 
Thus the proof of Lemma \ref{310} is complete. 
\end{proof}

In view of Lemma \ref{310}, we only need to consider the 
decoupled system \eqref{3a}, which makes the derivation of our 
energy estimates easier. 
As a priori estimates of solutions, we show the uniform energy 
estimate for $u$ by using the equation $\eqref{3a}_1$ and then 
a similar estimate for $q$ by making use of the relation 
$\eqref{3a}_2$. 
To this end, we consider solutions $u$ to the problem 
$\eqref{3a}_1$, \eqref{1a1}, which are defined on the time interval 
$[0,T]$ for $T>0$ and satisfy 
\be\label{3b}
\|u(t)\|_{L^{\infty}}\leq C,\qquad
\|\p_xu(t)\|_{L^{\infty}}\leq C
\ee 
for $0\leq t\leq T$. Notice that these estimates hold true if 
$\sup_{0\leq t\leq T}\|u(t)\|_{H^s}\leq C$ for $s\geq[{n\over2}]+2$. 

Now we multiply $\eqref{3a}_1$ by $u$. A direct computation gives 
\begin{equation*}
\begin{split}
&(u^2+|\nabla u|^2)_t+2|\nabla u|^2
   -2\nabla\cdot\{u\nabla(u_t+u+\sum_jf_j(u)_{x_j})\} \\
&+\sum_j\{2g_j(u)+f_j'(u)|\nabla u|^2\}_{x_j}
   =-\sum_jf_j''(u)u_{x_j}|\nabla u|^2,
\end{split}
\end{equation*}
where $g_j(u)=\int_0^uf_j'(\eta)\eta\,d\eta$. 
Integrating this equality with respect to $x$, we have 
\be\label{3c} 
{d\over{dt}}\|u(t)\|^2_{H^1}+2\|\partial_xu(t)\|^2_{L^2}
\leq C\|\p_xu(t)\|_{L^{\infty}}\|\partial_xu(t)\|^2_{L^2}. 
\ee
To get similar estimates for the derivatives, we apply 
$\partial_x^l$ to $\eqref{31}_1$, obtaining 
\begin{equation*}
\begin{split}
&\p_x^lu_t-\Delta\p_x^lu_t-\Delta\p_x^lu
   +\sum_jf_j'(u)\p_x^lu_{x_j}
   -\nabla\cdot\sum_jf_j'(u)\nabla\p_x^lu_{x_j} \\
&=-h^l+\nabla\cdot H^l,
\end{split}
\end{equation*}
where $h^l=\sum_j[\p_x^l,f_j'(u)]u_{x_j}$ and 
$H^l=\sum_j[\p_x^l\p_{x_j},f_j'(u)]\nabla u$. 
We multiply this equation by $\p_x^lu$ and compute directly to get 
\begin{equation*}
\begin{split}
&(|\p_x^lu|^2+|\nabla\p_x^lu|^2)_t+2|\nabla\p_x^lu|^2
   -2\nabla\cdot\{\p_x^lu\nabla\p_x^l(u_t+u) \\[1mm]
&\qquad\qquad
   +\p_x^lu\sum_jf_j'(u)\nabla\p_x^lu_{x_j}\}
   +\sum_j\{f_j'(u)(|\p_x^lu|^2+|\nabla\p_x^lu|^2)\}_{x_j} \\
&=\sum_jf_j''(u)u_{x_j}(|\p_x^lu|^2+|\nabla\p_x^lu|^2)
   -2(h^l\p_x^lu+H^l\cdot\nabla\p_x^lu)
   +2\nabla\cdot(H^l\p_x^lu).
\end{split}
\end{equation*}
Integrating this equality with respect $x$ and estimating the right 
hand side by applying Lemma \ref{31}, we obtain 
\be\label{3d}
{d\over{dt}}\|\p^l_xu(t)\|^2_{H^1}+2\|\p^l_x\nabla u(t)\|^2_{L^2}
\leq C\|\p_xu(t)\|_{L^{\infty}}\|\p^l_xu(t)\|^2_{H^1},
\ee
where $1\leq l\leq s-1$. Here we have used the fact that  $f_j(u)=O(u^2), j=1, \cdots, n$ for $u\to 0$. 
We add \eqref{3c} and \eqref{3d} for $1\leq l\leq s-1$ and 
integrate over $(0,t)$. This yields 
%
%
%
\be\label{3g-a}
\begin{split}
&\|u(t)\|^2_{H^{s}}+\int^t_0\|\p_xu(\tau)\|^2_{H^{s-1}}d\tau \\
&\leq CE_0^2
   +C\!\int_0^t\|\p_xu(\tau)\|_{L^{\infty}}
   \|\p_xu(\tau)\|^2_{H^{s-1}}d\tau,
\end{split}
\ee 
where $E_0=\|u_0\|_{H^s}$. 
Let $\bar{\delta}$ be a positive number (independnt of $T$) 
and assume that 
$\sup_{0\leq t\leq T}\|u(t)\|_{H^s}\leq\bar{\delta}$, where 
$s\geq[{n\over2}]+2$. Then the second term on the right hand side 
of \eqref{3g-a} is estimated by 
$C\bar{\delta}\int_0^t\|\p_xu(\tau)\|^2_{H^{s-1}}d\tau$. 
Therefore, choosing $\bar{\delta}$ so small that 
$C\bar{\delta}\leq {1\over2}$, we arrive at the uniform 
energy estimate 
\be\label{3g}
\|u(t)\|^2_{H^{s}}+\int^t_0\|\p_xu(\tau)\|^2_{H^{s-1}}d\tau 
\leq CE_0^2.
\ee 
On the other hand, it follows from $\eqref{3a}_2$ that 
$\|q\|_{H^{s+1}}\leq\|\p_xu\|_{H^{s-1}}$, which combined with 
\eqref{3g} yields 
\be\label{3h} 
\|q(t)\|^2_{H^{s+1}}+\int^t_0\|q(\tau)\|^2_{H^{s+1}}d\tau
\leq CE_0^2. 
\ee
These observations are summarized as follows. 

\bp\label{311}
Let $n\geq1$ and let $s\geq[{n\over2}]+2$ be an integer. 
Assume that $u_0\in H^{s}(\R^n)$ and put $E_0:=\|u_0\|_{H^{s}}$. 
Let $(u,q)$ be a solution to the problem \eqref{1a}, \eqref{1a1} 
on the time interval $[0,T]$. 
Then there is a small positive constant $\bar{\delta}$ independent 
of $T$ such that if 
$\sup_{0\leq t\leq T}\|u(t)\|_{H^s}\leq\bar{\delta}$, then the 
solution verifies the following uniform energy estimate for 
$t\in[0,T]$: 
\begin{equation}\label{3.9}
\|u(t)\|^2_{H^{s}}+\|q(t)\|^2_{H^{s+1}}
+\int^t_0\|\p_xu(\tau)\|^2_{H^{s-1}}+\|q(\tau)\|^2_{H^{s+1}}d\tau
\leq CE_0^2. 
\end{equation}
\ep

By virtue of the a priori estimate \eqref{3.9} for small solutions 
stated in Proposition \ref{311}, we can apply the continuity argument 
and obtain a unique global solution to the problem \eqref{1a}, 
\eqref{1a1}, provided that $E_0$ is suitably small, say, 
$E_0\leq\d_0$. The solution obtained verifies \eqref{3.9} for 
$t\geq 0$. This proves Theorem \ref{21}.


\subsection{Optimal decay estimates}

In this subsection, we obtain the optimal decay estimates of the 
solution constructed in Theorem \ref{21} by using the 
time-weighted energy method. 
To this end, we define two time-weighted energy norms $E(T)$ and $M(T)$. 
Also, we introduce $D(T)$ as the dissipation norm corresponding to 
$E(T)$. 
\begin{equation}\label{E}
E(T)^2:=\sum\limits^{{{s}}}\limits_{j=0}\sup\limits_{0\leq t\leq T}
(1+t)^{j}\|\p^{j}_xu(t)\|^2_{H^{s-j}}, 
\end{equation}
\begin{equation}\label{D}
D(T)^2:=\sum\limits^{{{s}}}\limits_{j=0}\int^T_0(1+\tau)^{j}
\|\p^{j+1}_xu(\tau)\|^2_{H^{s-j-1}}d\tau, 
\end{equation}
\begin{equation}\label{M}
M(T)^2:=\sum\limits^{{{s-1}}}\limits_{j=0}\sup\limits_{0\leq t\leq T}
(1+t)^{{n\over2}+j}\|\p^{j}_xu(t)\|^2_{H^{s-j}}, 
\end{equation}
where $s\geq1$. 
To derive estimates for $E(T)$, $D(T)$ and $M(T)$, we make use of 
the following time-weighted norm $N(T)$: 
\be\label{N}
N(T):=\sup\limits_{0\leq t\leq T}(1+t)\|\p_xu(t)\|_{L^{\infty}}.
\ee
As for the energy $E(T)$ and $D(T)$, we have the following estimate. 

\bp \label{32-3} 
Let $n\geq 1$ and $s\geq[{n\over2}]+2$. Assume that 
$u_0\in H^{s}(\R^n)$ and put $E_0:=\|u_0\|_{H^{s}}$.  
Then the solution to the problem \eqref{1a}, \eqref{1a1} 
constructed in Theorem \ref{21} satisfies following energy 
estimate: 
$$
E(T)^2+D(T)^2\leq CE_0^2+CN(T)D(T)^2.
$$ 
\ep

\begin{proof}
In order to prove this proposition, it is enough to show the following 
estimates for any $t\in[0,T]$ and $0\leq j\leq s$: 
\be\label{32-c}
\bal
&(1+t)^{j}\|\p^{j}_xu(t)\|^2_{H^{s-j}}
   +\int^t_0(1+\tau)^{j}\|\p^{j+1}_xu(\tau)\|^2_{H^{s-j-1}}d\tau \\[3mm]
&\leq CE_0^2+CN(T)D(T)^2.
\ea
\ee
We know from \eqref{2.1} that \eqref{32-c} holds true for $j=0$. 
Now, let $0\leq k\leq s-1$ and suppose that \eqref{32-c} holds true 
for $j=k$. Then we show \eqref{32-c} for $j=k+1$. 
Multiplying \eqref{3d} by $(1+t)^{k+1}$, integrating with respect 
to $t$ over $(0,t)$ and adding for $l$ with $k+1\leq l\leq s-1$, 
we have 
$$ 
\bal
&(1+t)^{k+1}\|\p_x^{k+1}u(t)\|^2_{H^{s-k-1}}
   +\int^t_0(1+\tau)^{k+1}\|\p_x^{k+2}u(\tau)\|^2_{H^{s-k-2}}d\tau 
\\[3mm]
&\leq CE_0^2
   +C\!\int^t_0(1+\tau)^{k}\|\p^{k+1}_xu(\tau)\|^2_{H^{s-k-1}}d\tau 
\\[3mm]
&\qquad
   +C\!\int^t_0(1+\tau)^{k+1}\|\p_xu(\tau)\|_{L^{\infty}}
   \|\p_x^{k+1}u(\tau)\|^2_{H^{s-k-1}}d\tau.
\ea
$$
The second term on the right hand side is estimated by the 
induction hypothesis \eqref{32-c} with $j=k$, while the last term 
can be estimated by $CN(T)D(T)^2$. Consequently, we have 
$$
\bal
&(1+t)^{k+1}\|\p^{k+1}_xu(t)\|^2_{H^{s-k-1}}
   +\int^t_0(1+\tau)^{k+1}
   \|\p_x^{k+2}u(\tau)\|^2_{H^{s-k-2}}d\tau \\[3mm]
&\leq CE_0^2+CN(T)D(T)^2,
\ea
$$
which shows that \eqref{32-c} holds true for $j=k+1$.
Thus, by induction, we have proved Proposition \ref{32-3}. 
\end{proof}

By employing the optimal decay results expressed in $E(T)$ and 
$M(T)$, we obtain the following estimates for $N(T)$. 

\bl\label{324}
{\rm (i)}\ If $n\geq2$ and $s\geq [{n\over2}]+2$, 
then $N(T)\leq CE(T)$. \\
\noindent
{\rm (ii)}\ If $n=1$ and $s\geq3$, then $N(T)\leq CM(T)$.
\el

\begin{proof}
(i)\ Let $s_0=[{n\over2}]+1$ and $\theta={n\over2s_0}$. 
By applying the Gagliardo-Nirenberg inequality, we see that 
$$
\bal
\|\p_xu(t)\|_{L^{\infty}}
&\leq C\|\p_xu(t)\|_{L^2}^{1-\theta}
   \|\p_x^{s_0+1}u(t)\|_{L^2}^{\theta} \\[2mm]
&\leq CE(T)(1+t)^{-{1\over2}(1-\theta)-{s_0+1\over2}\theta} \\[2mm]
&=CE(T)(1+t)^{-({n\over4}+{1\over2})}.
\ea
$$ 
Since $n\geq2$, it yields that $N(T)\leq CE(T)$. 

(ii)\ By using the one-dimensional Gagliardo-Nirenberg inequality, 
we have 
\begin{equation*} 
\|\p_xu(t)\|_{L^{\infty}}
\leq C\|\p_xu(t)\|_{L^2}^{{1\over2}}
   \|\p^{2}_xu(t)\|^{{1\over2}}_{L^2}\leq CM(T)(1+t)^{-1},
\end{equation*}
which gives $N(T)\leq CM(T)$. This completes the proof. 
\end{proof}

\begin{proof}[Proof of Corollary \ref{21a}]
By virtue of Lemma \ref{324} (i), we have $N(T)\leq CE(T)$ for 
$n\geq 2$, which together with Proposition \ref{32-3} gives 
$$
E(T)^2+D(T)^2\leq CE_0^2+CE(T)D(T)^2. 
$$
Put $X(T):=E(T)+D(T)$. Then we have $X(T)^2\leq CE_0^2+CX(T)^3$. 
This inequality is solved as $X(T)\leq CE_0$, provided that 
$E_0$ is suitably small. In particular, we have $E(T)\leq CE_0$, 
which shows the desired decay estimate \eqref{2.2} for 
$0\leq k\leq s$. 
Moreover, by virtue of $\eqref{3a}_2$, we have 
$$
\|\p^k_xq(t)\|_{H^{s+1-k}}
\leq\|\p_x^{k+1}u(t)\|_{H^{s-1-k}}\leq CE_0(1+t)^{-{k+1\over2}},
$$ 
for $0\leq k\leq s-1$, where we have used \eqref{2.2} with 
$k$ replaced by $k+1$. This completes the proof of Corollary 
\ref{21a}. 
\end{proof}

To estimate the energy $M(T)$, we need the following $L^1$ estimate 
of the solution. 

\bl\label{320}
Under the same assumptions as in Theorem \ref{21b}, 
the solution to the problem \eqref{1a}, \eqref{1a1} satisfies 
the following $L^1$ estimate for $u$: 
\begin{equation}\label{L^1}
\|u(t)\|_{L^1}\leq\|u_0\|_{L^1}.
\end{equation}
\el

\begin{proof}
Applying $(1-\Delta)^{-1}$ to $\eqref{3a}_1$, we have 
\begin{equation}\label{3aa}
u_t+\sum\limits_{j=1}\limits^{n}f^j(u)_{x_j}+u-(I-\Delta)^{-1} u=0.
\end{equation}
We denote by $K(x)$ the the fundamental solution to the operator 
$I-\Delta$, that is, $(I-\Delta)^{-1} u=K\ast u$. We know that 
$K(x)\geq 0$, $K\in L^1$ and $\int_{\R^n}K(x)dx=1$. 
See \cite{Di} for the details. 
Let $j_{\delta}$ be the Friedrichs mollifier and put 
$$
\phi_{\delta}(u)=j_{\delta}\ast {\rm sign}(u), \quad
\Phi_{\delta}(u)=\int^u_0\phi_{\delta}(\xi)d\xi.
$$
We multiply \eqref{3aa} by $\phi_{\delta}(u)$ and integrate the 
resulting equation over $\R^n\times(0,t)$. Then, letting 
$\delta\rightarrow 0$, we obtain the desired $L^1$ estimate 
\eqref{L^1} just in the same way as in \cite{Di,KT,RZ}. 
The details are omitted. 
\end{proof}

By employing the time-weighted energy method together with 
the $L^1$ estimate \eqref{L^1}, we get the following estimate 
for $M(T)$. 

\bp \label{323} 
Let $n\geq 1$ and $s\geq[{n\over2}]+2$. Assume that 
$u_0\in H^{s}(\R^n)\cap L^1(\R^n)$, and put $E_0:=\|u_0\|_{H^s}$ and 
$E_1:=\|u_0\|_{H^{s}}+\|u_0\|_{L^1}$. Then, if $E_0$ is suitably 
small, then the solution to the problem \eqref{1a}, \eqref{1a1} 
constructed in Theorem \ref{21} satisfies following estimate: 
$$
M(T)^2\leq CE_1^2(1+ N(T))^{s-1}.
$$ 
\ep

\begin{proof}
In order to prove this proposition, it is enough to show the 
following estimate for any $t\in[0,T]$ and $0\leq j\leq s-1$: 
\be\label{32b} 
\|\p^{j}_xu(t)\|^2_{H^{s-j}} \leq
CE_1^2(1+N(T))^{s-1}(1+t)^{-{n\over2}-j}.
\ee 
For this purpose, it is sufficient to prove 
\be\label{32c}
\bal
&(1+t)^{j+\a}\|\p^{j}_xu(t)\|^2_{H^{s-j}}
   +\int^t_0(1+\tau)^{j+\a}
   \|\p^{j+1}_xu(\tau)\|^2_{H^{s-j-1}}d\tau \\[3mm]
&\leq CE_1^2(1+N(T))^j(1+t)^{\a-{n\over2}}
\ea
\ee
for $t\in[0,T]$ and $0\leq j\leq s-1$, where $\a>{n\over2}$. 

First we show \eqref{32c} for $j=0$. 
We add \eqref{3c} and \eqref{3d} for $1\leq l\leq s-1$, multiply 
the resulting inequality by $(1+t)^{\a}$, and then integrate over 
$(0,t)$. This yields 
\begin{equation}
\bal
&(1+t)^{\a}\|u(t)\|^2_{H^{s}}
   +\int^t_0(1+\tau)^{\a}\|\p_xu(\tau)\|^2_{H^{s-1}}d\tau \\[3mm]
&\leq CE_0^2+C\!\int^t_0
   (1+\tau)^{\a-1}\|u(\tau)\|^2_{H^{s}}d\tau \\[3mm]
&\quad
   +C\!\int^t_0(1+\tau)^{\a}\|\p_xu(\tau)\|_{L^{\infty}}
   \|\p_xu(\tau)\|^2_{H^{s-1}}d\tau 
=:CE_0^2+I_1+I_2.
\ea
\end{equation}
Since $\|u(t)\|_{H^s}\leq CE_0$ by \eqref{2.1}, we can estimate 
the term $I_2$ as 
$$
\bal
I_2
\leq CE_0\!\int^{t}_{0}(1+\tau)^{\a}\|\p_xu(\tau)\|^2_{H^{s-1}}d\tau
\leq {1\over4}\int^{t}_{0}(1+\tau)^{\a}\|\p_xu(\tau)\|^2_{H^{s-1}}d\tau,
\ea
$$
where $E_0$ is assumed to be small as $CE_0\leq{1\over4}$. 
On the other hand, we divide $I_1$ into two parts: 
$$
\bal
I_1=C\!\int^t_0(1+\tau)^{\a-1}\|u(\tau)\|^2_{L^{2}}d\tau
   +C\!\int^t_0(1+\tau)^{\a-1}\|\p_xu(\tau)\|^2_{H^{s-1}}d\tau
=:I_{11}+I_{12}.
\ea
$$
To estimate $I_{12}$, we choose $T_1$ so large that 
$C(1+T_1)^{-1}\leq{1\over4}$. Then we divide the time interval 
$[0,t]$ into two parts $[0,T_1]$ and $[T_1,t]$; here we treat the 
case $t\geq T_1$ because the case $t\leq T_1$ is easier. 
Thus we have 
$$
\bal
I_{12}
&\leq C\max\limits_{0\leq t\leq T_1}\{(1+t)^{\a-1}\}
   \!\int^{T_1}_0\|\p_xu(\tau)\|^2_{H^{s-1}}d\tau
	   +{1\over4}\int^{t}_{T_1}(1+\tau)^{\a}
   \|\p_xu(\tau)\|^2_{H^{s-1}}d\tau \\[3mm]
&\leq CE_0^2
   +{1\over4}\int^{t}_{0}(1+\tau)^{\a}
   \|\p_xu(\tau)\|^2_{H^{s-1}}d\tau,
\ea
$$
where we have used \eqref{2.1}. 
Finally, we estimate the term $I_{11}$. 
By using the Gagliardo-Nirenberg inequality 
$\|u\|_{L^2}\leq C\|\p_xu\|_{L^2}^{\theta}\|u\|_{L^1}^{1-\theta}$ 
with $\theta={n\over n+2}$ and applying the Young inequality, 
we can estimate $I_{11}$ as 
$$
\bal
I_{11}
&\leq C\!\int_0^t(1+\tau)^{\a-1}
   \|\p_xu(\tau)\|_{L^2}^{2\theta}
	   \|u(\tau)\|_{L^1}^{2(1-\theta)}d\tau \\[3mm]
&\leq {1\over4}\int^t_0(1+\tau)^{\a}\|\p_xu(\tau)\|^2_{L^{2}}d\tau
   +C\!\int_0^t(1+\tau)^{\a-{n\over2}-1}\|u(\tau)\|_{L^1}^2d\tau \\[3mm]
&\leq {1\over4}\int^t_0(1+\tau)^{\a}\|\p_xu(\tau)\|^2_{H^{s-1}}d\tau
   +C\|u_0\|_{L^1}^2(1+t)^{\a-{n\over2}}, 
\ea
$$
where we have used the $L^1$ estimate \eqref{L^1} and the 
condition $\a>{n\over2}$. 
Consequently, under the smallness assumption on $E_0$, we arrive at 
the estimate 
$$
(1+t)^{\a}\|u(t)\|^2_{H^{s}}
   +\int^t_0(1+\tau)^{\a}\|\p_xu(\tau)\|^2_{H^{s-1}}d\tau
\leq CE_1^2(1+t)^{\a-{n\over2}},
$$ 
which proves \eqref{32c} for $j=0$. 

Now, let $0\leq k\leq s-2$ and suppose that \eqref{32c} holds true 
for $j=k$. Then we show \eqref{32c} for $j=k+1$. 
Multiplying \eqref{3d} by $(1+t)^{k+1+\a}$, integrating with respect 
to $t$ over $(0,t)$ and adding up for $l$ with $k+1\leq l\leq s-1$, 
we have 
$$ 
\bal
&(1+t)^{k+1+\a}\|\p_x^{k+1}u(t)\|^2_{H^{s-k-1}}
   +\int^t_0(1+\tau)^{k+1+\a}
   \|\p_x^{k+2}u(\tau)\|^2_{H^{s-k-2}}d\tau \\[3mm]
&\leq CE_0^2+C\!\int^t_0(1+\tau)^{k+\a}
   \|\p^{k+1}_xu(\tau)\|^2_{H^{s-k-1}}d\tau \\[3mm]
&\qquad
   +C\!\int^t_0(1+\tau)^{k+1+\a}
   \|\p_xu(\tau)\|_{L^{\infty}}\|\p_x^{k+1}u(\tau)\|^2_{H^{s-k-1}}d\tau.
\ea
$$
Here, using the induction hypothesis \eqref{32c} with $j=k$, we can 
estimate the second term on the right hand side by 
$CE_1^2(1+N(T))^k(1+t)^{\a-{n\over2}}$. 
Similarly, we estimate the last term as 
$$ 
\bal
&C\!\int^t_0(1+\tau)^{k+1+\a}
   \|\p_xu(\tau)\|_{L^{\infty}}
   \|\p_x^{k+1}u(\tau)\|^2_{H^{s-k-1}}d\tau \\[3mm]
&\leq CN(T)\!\int^t_0(1+\tau)^{k+\a}
   \|\p_x^{k+1}u(\tau)\|^2_{H^{s-k-1}}d\tau \\[3mm]
&\leq CE_1^2N(T)(1+N(T))^k(1+t)^{\a-{n\over2}}.
\ea
$$
Thus we obtain 
$$
\bal
&(1+t)^{k+1+\a}\|\p^{k+1}_xu(t)\|^2_{H^{s-k-1}}
   +\int^t_0(1+\tau)^{k+1+\a}
   \|\p_x^{k+2}u(\tau)\|^2_{H^{s-k-2}}d\tau \\[3mm]
&\leq CE_1^2(1+N(T))^{k+1}(1+t)^{\a-{n\over2}}.
\ea
$$
This shows that \eqref{32c} holds true for $j=k+1$. 
Thus, by induction, we have proved Proposition \ref{323}. 
\end{proof}

\begin{proof}[Proof of Theorem \ref{21b}]
By virtue of Lemma \ref{324}, we have $N(T)\leq C(E(T)+M(T))$. 
Therefore we have from Propositions \ref{32-3} and \ref{323} that 
$$
\bal
&E(T)^2+D(T)^2\leq CE_0^2+C(E(T)+M(T))D(T)^2, \\[2mm]
&M(T)^2\leq CE_1^2(1+E(T)+M(T))^{s-1}. 
\ea
$$
Put $Y(T):=E(T)+D(T)+M(T)$. Then we have 
$Y(T)^2\leq CE_1^2(1+Y(T))^{s-1}+CY(T)^3$. 
This inequality is solved as $Y(T)\leq CE_1$, provided that 
$E_1$ is suitably small, say, $E_1\leq\delta_1$. 
In particular, we have $M(T)\leq CE_1$, which proves the decay 
estimate \eqref{2.4} for $0\leq k\leq s-1$. 
Moreover, using $\eqref{3a}_2$ and \eqref{2.4} with $k$ replaced by
$k+1$, we obtain 
$$
\|\p^k_xq(t)\|_{H^{s+1-k}}
\leq\|\p_x^{k+1}u(t)\|_{H^{s-1-k}}
\leq CE_1(1+t)^{-{{n\over4}+{k+1\over2}}}
$$ 
for $0\leq k\leq s-2$. This completes the proof of Theorem 
\ref{21b}. 
\end{proof}


\section{Asymptotic profile}

The aim of this section is to prove Theorem \ref{22} on the 
asymptotic profile. To this end, we first consider the 
corresponding linear problem \eqref{32a}. 
The fundamental solution to $\eqref{32a}_1$ is given by 
$G(x,t)=\mathcal{F}^{-1}[e^{-{|\xi|^2\over 1+|\xi|^2}t}](x)$, and 
the solution to $\eqref{32a}_1$ can be expressed in terms of the 
fundamental solution as 
\begin{equation}\label{LS}
\bar{u}(x,t)=(G(t)\ast u_0)(x).
\end{equation}

The solution operator $G(t)*$ verifies the following decay 
property: 

\bl\label{41}
Let $n\geq 1$ and $s\geq 0$. 
If $\phi\in H^s(\R^n)\cap L^1(\R^n)$, then we have the following 
decay estimate: 
$$ 
\|\p_x^kG(t)\ast \phi\|_{L^2}
\leq C(1+t)^{-{n\over4}-{k\over2}}\|\phi\|_{L^1}
+Ce^{-{t\over2}}\|\p_x^k\phi\|_{L^2}, \quad
0\leq k\leq s.
$$
\el

\begin{proof}
By direct calculation, we have 
$$ 
\bal
\|\p_x^kG(t)\ast\phi\|^2_{L^2} 
&\leq C\!\left(\int_{|\xi|\leq1}+\int_{|\xi|>1}\right)
   |\xi|^{2k}e^{-{2|\xi|^2\over1+|\xi|^2}t}|\hat{\phi}(\xi)|^2d\xi \\[3mm]
&\leq C\!\int_{|\xi|\leq1}|\xi|^{2k}
   e^{-{|\xi|^2}t}|\hat{\phi}(\xi)|^2d\xi
   +C\int_{|\xi|>1}|\xi|^{2k}e^{-t}|\hat{\phi}(\xi)|^2d\xi \\[3mm]
&\leq C\!\int_{|\xi|\leq1}|\xi|^{2k}
   e^{-{|\xi|^2}(t+1)}e^{|\xi|^2}|\hat{\phi}(\xi)|^2d\xi
   +C\int_{|\xi|>1}|\xi|^{2k}e^{-t}|\hat{\phi}(\xi)|^2d\xi \\[3mm]
&\leq C(1+t)^{-{n\over2}-k}\|\phi\|^2_{L^1}
   +Ce^{-t}\|\p_x^k\phi\|^2_{L^2}.
\ea
$$ 
This completes the proof. 
\end{proof}

By Duhamel principle, we can express the solution of $\eqref{3a}_1$ 
(or \eqref{3aa}), \eqref{1a1} as follows: 
\be\label{4.2}
u(x,t)=G(t)\ast u_0
-\int^t_0G(t-\tau)\ast{\rm div}f(u)(\tau)d\tau.
\ee
Here and in the following, we use the abbreviation 
$f(u)=(f_1(u),\cdots, f_n(u))$.
We decompose the solution formula \eqref{4.2} in the form 
$u(t)=\bar{u}(t)-F(u)(t)$, where $\bar{u}$ is the linear solution 
given in \eqref{LS} and 
\begin{equation}\label{4.3}
F(u)(x,t)=\int^t_0G(t-\tau)\ast{\rm div}f(u)(\tau)d\tau.
\end{equation}

First we prove that the solution to the problem 
$\eqref{3a}_1$, \eqref{1a1} can be approximated by the solution to 
the corresponding linear problem $\eqref{32a}_1$. 

\bp\label{NL}
Let $n\geq 2$ and $s\geq [{n\over2}]+2$. Assume that 
$u_0\in H^{s}(\R^n)\cap L^1(\R^n)$ 
and put $E_1:=\|u_0\|_{H^{s}}+\|u_0\|_{L^1}$. 
Let $(u,q)$ be the global solution to the problem \eqref{1a}, 
\eqref{1a1} which is obtained in Theorem \ref{21}, and let 
$\bar{u}$ be the solution to the corresponding linear problem \eqref{32a},
which is given by the formula \eqref{LS}. Then we have 
$$
\|\p^k_x(u-\bar{u})(t)\|_{H^{s-1-k}}
\leq CE_1^2\rho(t)(1+t)^{-{n\over4}-{k+1\over2}}
$$ 
for $k$ with $0\leq k\leq s-1$, where 
$\rho(t)=\ln(1+t)$ for $n=2$ and $\rho(t)=1$ for $n\geq3$. 
\ep

\begin{proof}
Let $k$ and $m$ be nonnegative integers. 
We apply $\p^{k+m}_x$ to $F(u)$ in \eqref{4.3} and take the $L^2$ 
norm, obtaining 
\begin{equation}\label{I}
\bal
\|\p_x^{k+m}F(u)(t)\|_{L^2}&\leq
\left(\int_0^{{t\over2}}+\int^t_{{t\over2}}\right)
   \|\p_x^{k+m+1}G(t-\tau)\ast f(u)(\tau)\|_{L^2}d\tau \\[3mm]
&=:I_1+I_2.
\ea
\end{equation}
By applying Lemma \ref{41}, using Lemma \ref{31} and noticing that $f_j(u)=O(u^2), j=1, \cdots, n$ for $u\to 0$, we have 
\begin{equation*}\label{I1}
\bal 
I_1&\leq
C\!\int^{{t\over2}}_0(1+t-\tau)^{-{n\over4}-{k+m+1\over2}}
   \|f(u)\|_{L^1}d\tau \\[2mm]
&\qquad\qquad
   +C\!\int^{{t\over2}}_0e^{-{t-\tau\over2}}
   \|\p_x^{k+m+1}f(u)\|_{L^2}d\tau \\[3mm]
&\leq
C\!\int^{{t\over2}}_0(1+t-\tau)^{-{n\over4}-{k+m+1\over2}}
   \|u\|^2_{L^2}d\tau \\[2mm]
&\qquad\qquad
   +C\!\int^{{t\over2}}_0e^{-{t-\tau\over2}}
   \|u\|_{L^{\infty}}\|\p_x^{k+m+1}u\|_{L^2}d\tau.
\ea
\end{equation*}
Similarly, we have 
\begin{equation*}\label{I2}
\bal 
I_2&\leq
C\!\int^t_{{t\over2}}(1+t-\tau)^{-{n\over4}-{1\over2}}
   \|\p^{k+m}_xf(u)\|_{L^1}d\tau \\[2mm]
&\qquad\qquad
   +C\!\int^t_{{t\over2}}e^{-{t-\tau\over2}}
   \|\p_x^{k+m+1}f(u)\|_{L^2}d\tau \\[3mm]
&\leq
C\!\int^t_{{t\over2}}(1+t-\tau)^{-{n\over4}-{1\over2}}
   \|u\|_{L^2}\|\p^{k+m}_xu\|_{L^2}d\tau \\[2mm]
&\qquad\qquad
   +C\!\int^t_{{t\over2}}e^{-{t-\tau\over2}}
   \|u\|_{L^{\infty}}\|\p_x^{k+m+1}u\|_{L^2}d\tau.
\ea
\end{equation*}

We estimate the terms $I_1$ and $I_2$. 
Let $s_0=[{n\over2}]+1$ and $\theta={n\over2s_0}$. 
By using the Gagliardo-Nirenberg inequality and \eqref{2.2}, 
we see that 
\be\label{B}
\bal 
\|u(t)\|_{L^{\infty}}&\leq
C\|u(t)\|_{L^2}^{1-\theta}\|\p^{s_0}_xu(t)\|^{\theta}_{L^2}\\[2mm]
&\leq
CE_1(1+t)^{-{s_0\over2}\theta}=CE_1(1+t)^{-{n\over4}}.
\ea
\ee
By using \eqref{2.2}, \eqref{2.4} and \eqref{B}, we estimate 
$I_1$ as 
\begin{equation}\label{I11}
\bal
I_1&\leq
CE_1^2\int^{{t\over2}}_0(1+t-\tau)^{-{n\over4}-{k+m+1\over2}}
   (1+\tau)^{-{n\over2}}d\tau \\[2mm]
&\qquad
   +CE_1^2\int^{{t\over2}}_0e^{-{t-\tau\over2}}
   (1+\tau)^{-{n\over4}-{k+1\over2}}d\tau \\[2mm]
&\leq CE_1^2\rho(t)(1+t)^{-{n\over4}-{k+1\over2}}
\ea
\end{equation} 
for $k$ satisfying $k+m+1\leq s$, where $\rho(t)$ is given above. 
Similarly, we can estimate $I_2$ as 
\begin{equation}\label{I22}
\bal
I_2&\leq
CE_1^2\int^t_{{t\over2}}(1+t-\tau)^{-{n\over4}-{1\over2}}
   (1+\tau)^{-{n\over2}-{k+m\over2}}d\tau \\[2mm]
&\qquad\qquad
   +CE_1^2\int^t_{{t\over2}}e^{-{t-\tau\over2}}
   (1+\tau)^{-{n\over4}-{k+1\over2}}d\tau \\[2mm]
&\leq CE_1^2\rho(t)(1+t)^{-{n\over4}-{k+1\over2}}
\ea
\end{equation}
for $k$ satisfying $k+m+1\leq s$. 

We substitute \eqref{I11} and \eqref{I22} into \eqref{I} and take 
the sum for $m$ with $0\leq m\leq s-k-1$. This yields 
$$
\|\p_x^{k}F(u)(t)\|_{H^{s-1-k}}\leq
 CE_1^2\rho(t)(1+t)^{-{n\over4}-{k+1\over2}}
$$ 
for $k$ with $0\leq k\leq s-1$. 
This completes the proof of Proposition \ref{NL}. 
\end{proof}

Next we consider the  simpler problem for the linear 
heat equation \eqref{HE}. 
Notice that the solution to $\eqref{HE}_1$ is expressed as 
$\tilde{u}(x,t)=(G_0(t)\ast u_0)(x)$, where $G_0(x,t)$ is 
the heat kernel given in \eqref{FS}. 
By direct computation we see that the solution $\bar{u}(x,t)$ 
of $\eqref{32a}_1$ is well approximated by the solution 
$\tilde{u}(x,t)$ of $\eqref{HE}_1$ as $t\to+\infty$:

\bl\label{LH}
Let $n\geq1$ and $s\geq 0$. Assume that 
$u_0\in H^{s}(\R^n)\cap L^1(\R^n)$ and put 
$E_1=\|u_0\|_{H^{s}}+\|u_0\|_{L^1}$. 
Let $\bar{u}(x,t)$ and $\tilde{u}(x,t)$ be solutions to 
$\eqref{32a}_1$ and $\eqref{HE}_1$, respectively. 
Then we have the following asymptotic relation: 
\begin{equation*}
\|\p^k_x(\bar{u}-\tilde{u})(t)\|_{H^{s-k}}
\leq CE_1(1+t)^{-{n\over4}-{k\over2}-{1}}, \quad 0\leq k\leq s.
\end{equation*}
\el
\begin{proof}
\noindent
Since 
$
|(\hat{G}(\xi,t)-\hat{G}_0(\xi,t)|\leq C|\xi|^2e^{-c|\xi|^2t}$ for $|\xi|\leq 1$, and 
 $|(\hat{G}(\xi,t)-\hat{G}_0(\xi,t)|\leq Ce^{-ct}$ for $|\xi|\geq 1$,
by direct calculation we have  the following decay property for $(G-G_0)(t)*$: $\forall \phi\in H^s\cap L^1$,
\begin{equation*}\label{}
\|\p^k_x(G-G_0)(t)\ast\phi\|_{L^2}
  \leq C(1+t)^{-{n\over4}-{k\over2}-1}\|\phi\|_{L^1}
  +Ce^{-ct}\|\p^k_x\phi\|_{L^2},
\end{equation*}
here $k\leq s$. Replace $k$ with $k+m$, $\phi$ with $u_0$, then we have 
\begin{equation*}\label{}
\bal
\|\p^{k+m}_x(G-G_0)(t)\ast u_0\|_{L^2}
  &\leq C(1+t)^{-{n\over4}-{k+m\over2}-1}\|u_0\|_{L^1}  +Ce^{-ct}\|\p^{k+m}_xu_0\|_{L^2}
\\&\leq C(1+t)^{-{n\over4}-{k\over2}-1}\|u_0\|_{L^1}  +Ce^{-ct}\|u_0\|_{H^{k+m}}.
\ea
\end{equation*}
Take sum with $0\leq m\leq s-k$, then we complete the proof of Lemma \ref{LH}.
\end{proof}

Finally, we show that the solution to \eqref{HE} is well 
approximated by the linear diffusion wave defined in \eqref{LDW}. 
For this purpose, we recall the following well known result 
for the heat kernel. 

\bl\label{EO} 
Let $n\geq 1$ and $1\leq q\leq 2$. If $\phi\in L^q(\R^n)$, then 
we have 
\begin{equation}\label{EO1}
\|\p^k_xG_0(t)*\phi\|_{L^2}
\leq Ct^{-{n\over2}({1\over q}-{1\over2})-{k\over2}}
\|\phi\|_{L^q}, \quad k\geq 0.
\end{equation}
Also, if $\phi\in L^1_1(\R^n)$ and $\int_{\mathbb{R}^n}\phi(x)dx=0$, 
then we have 
\begin{equation}\label{EO2}
\|\p^k_xG_0(t)*\phi\|_{L^2}
  \leq Ct^{-{n\over4}-{{k+1}\over2}}\|\phi\|_{L^1_1}, \quad k\geq 0.
\end{equation}
\el


\begin{proof}[Proof of Theorem \ref{22}]
We recall the linear diffusion wave $(u^*,q^*)(x,t)$ in \eqref{LDW}. 
Since $G_0(x,t+1)$ is a solution of the linear heat equation 
$u_t-\Delta u=0$ with the initial data 
$G_0(x,1)=(4\pi)^{-{n\over2}}e^{-{|x|^2\over4}}=:\phi_0(x)$, 
we see that $G_0(x,t+1)=(G_0(t)*\phi_0)(x)$. 
Consequently, we can write 
$$
(\tilde{u}-u^*)(t)=G_0(t)*(u_0-M\phi_0),
$$
where $M=\int_{\mathbb{R}^n}u_0(x)dx$. 
Here $\phi_0(x)$ is a rapidly decreasing function satisfying 
$\int_{\mathbb{R}^n}\phi_0(x)dx=1$, so that we have 
$\int_{\mathbb{R}^n}(u_0-M\phi_0)(x)dx=0$. 
Therefore, applying Lemma \ref{EO}, we deduce that 
\begin{equation}\label{W1}
\begin{split}
\|\p_x^k(\tilde{u}-u^*)(t)\|_{H^{s-k}}
&=\|\p^k_xG_0(t)*(u_0-M\phi_0)\|_{H^{s-k}} \\[1mm]
&\leq CE_2(1+t)^{-{n\over4}-{{k+1}\over2}}
\end{split}
\end{equation}
for $0\leq k\leq s$. Now we write 
$u-u^*=(u-\bar{u})+(\bar{u}-\tilde{u})+(\tilde{u}-u^*)$ 
and apply Proposition \ref{NL}, Lemma \ref{LH} and \eqref{W1}. 
This yields the desired estimate \eqref{asymp-1} for 
$0\leq k\leq s-1$. 

On the other hand, we have $q^*=-\nabla u^*$ by \eqref{LDW}, 
which can be rewritten as 
$q^*=-(1-\Delta)^{-1}\nabla u^*+(1-\Delta)^{-1}\Delta\nabla u^*$.
Therefore, noting \eqref{3a},  we 
find that 
$$
q-q^*=-(1-\Delta)^{-1}\nabla(u-u^*)
-(1-\Delta)^{-1}\Delta\nabla u^*.
$$
Consequently, we have 
\begin{equation*}
\|\p_x^k(q-q^*)(t)\|_{H^{s-k}}
\leq \|\p_x^{k+1}(u-u^*)(t)\|_{H^{s-2-k}}
   +\|\p_x^{k+3}u^*(t)\|_{H^{s-2-k}}. 
\end{equation*}
Here, applying \eqref{asymp-1} with $k$ replaced by $k+1$, we see 
that the first term on the right hand side is estimated by 
$CE_2\rho(t)(1+t)^{-{n\over4}-{{k+2}\over2}}$ for 
$0\leq k\leq s-2$, where $\rho(t)$ is given in Theorem \ref{22}. 
Also, using the expression $u^*=MG_0(t)*\phi_0$ and 
applying \eqref{EO1}, we know that the second term can be 
estimated by $C\|u_0\|_{L^1}(1+t)^{-{n\over4}-{{k+3}\over2}}$. 
These observations prove the estimate \eqref{asymp-2}. 
Thus the proof of Theorem \ref{22} is complete. 
\end{proof}


\noindent{\bf Acknowledgments.}\ 
This work was partially supported by Grant-in-Aid for JSPS Fellows.


\end{document}